\documentclass[11pt]{article}
\usepackage{epsfig}
\usepackage{t1enc}
    \usepackage[latin1]{inputenc}
    \usepackage[english]{babel}
        \usepackage{latexsym}
\usepackage{amssymb}
\usepackage{euscript}
\begin{document}

\setlength{\arraycolsep}{.136889em}
\renewcommand{\theequation}{\thesection.\arabic{equation}}
\newtheorem{thm}{Theorem}[section]
\newtheorem{propo}{Proposition}[section]
\newtheorem{lemma}{Lemma}[section]
\newtheorem{corollary}{Corollary}[section]
\newtheorem{remark}{Remark}[section]
\def\begg{\begin{equation}}
\def\endd{\end{equation}}
\def\ep{\varepsilon}
\def\noo{n\to\infty}
\def\al{\alpha}
\def\be{\bf E}
\def\bp{\bf P}
\medskip
\centerline{\Large\bf Random walks on comb-type subsets of $\mathbb
Z^2$}

\bigskip\bigskip

\bigskip\bigskip

\renewcommand{\thefootnote}{1}

\noindent 
{\textbf{Endre Cs\'aki}}

\noindent
Alfr\'ed R\'enyi Institute of Mathematics, Hungarian Academy of
Sciences, Budapest, P.O.B. 127, H-1364, Hungary. E-mail address:
csaki.endre@renyi.mta.hu

\bigskip
\renewcommand{\thefootnote}{1}
\noindent 
{\textbf{Ant\'onia F\"oldes}\footnote{Research
supported by a PSC CUNY Grant, No. 61520-0049}}

\noindent
Department of Mathematics, College of Staten Island, CUNY, 2800
Victory Blvd., Staten Island, New York 10314, U.S.A.  E-mail
address: antonia.foldes@csi.cuny.edu

\bigskip

\centerline{\bf Abstract}

\medskip
\noindent We study the path behavior of the
simple symmetric  walk  on some comb-type subsets of ${\mathbb Z}^2$
which are obtained from ${\mathbb Z}^2$ by removing all horizontal
edges belonging to certain sets of values on the $y$-axis. We obtain some 
strong approximation results and  discuss their consequences.

\medskip


\medskip

 \vspace{.1cm}

\section{Introduction}

\renewcommand{\thesection}{\arabic{section}} \setcounter{equation}{0}

\setcounter{thm}{0} \setcounter{lemma}{0}

 An anisotropic walk is defined  as a nearest neighbor random walk 
on the square lattice $\mathbb{Z}^2$ of the plane with possibly unequal 
symmetric horizontal and vertical step probabilities, so that these 
probabilities depend only on the value of the vertical coordinate. 
More formally, consider the random 
walk $\{{\bf C}(N)=\left(C_1(N),C_2(N)\right);\, N=0,1,2,\ldots\}$
on $\mathbb{Z}^2$ with the transition probabilities
\begin{eqnarray} 
\label{eqn:2}
\nonumber
{\bf P}({\bf C}(N+1)&=&(k+1,j)|{\bf C}(N)=(k,j))={\bf P}
({\bf C}(N+1)=(k-1,j)|{\bf C}(N)=(k,j))=\frac{1}{2}-p_j, \\ 
 {\bf P} ({\bf C}(N+1)&=&(k,j+1)|{\bf C}(N)=(k,j))={\bf P} 
({\bf C}(N+1)=(k,j-1)|{\bf C}(N)=(k,j))=p_j, \label{prob1}
\end{eqnarray}
for $(k,j)\in{\mathbb Z^2}$, $N=0,1,2,\ldots$ with $0<p_j\leq 1/2$ and 
$\min_{j\in\mathbb Z} p_j<1/2$. 
Unless otherwise stated, we assume also that 
${\mathbf C}(0)=(0,0)$. 

\bigskip
In the present paper we are  interested in a special type of this anisotropic 
walk. We only want 
to consider walks for which  $p_j$ in (\ref{prob1})  is either $1/2$ or $1/4.$  
In particular, for such walks we consider an arbitrary subset $B$  of the 
integers on the $y$-axis and remove from the two-dimensional integer lattice 
all the horizontal lines which do not belong to the $y$-levels in $B.$ Denote 
this lattice by ${\mathbb C}^2={\mathbb C}^2(B)$. The transition probabilities 
throughout this paper are

\begin{eqnarray}
\nonumber
p_y&=&{\bf P}({\bf C}(N+1)=(x\pm 1,y)\mid {\bf C}(N)=(x,y))\\ 
\nonumber
&= &{\bf P}({\bf C}(N+1)=(x,y\pm 1)\mid {\bf C}(N)=(x,y))=\frac{1}{4},
\quad {\rm if} \, \, y\in B\\
p_y&=&{\bf P}({\bf C}(N+1)=(x,y\pm 1)\mid {\bf C}(N)=(x,y))=\frac12,
\quad {\rm if} \,\,  y   \notin B , \label{combl}
\end{eqnarray}

A compact way of describing the just introduced transition
probabilities for this simple random walk ${\bf C}(N)$ on ${\mathbb
C}^2(B)$ is via defining

\begin{equation}
p({\bf u,v}):={\bf P}({\bf C}(N+1)={\bf v}\mid {\bf C}(N)={\bf u})=
\frac1{{\rm deg}({\bf u})},
\end{equation}
for locations ${\bf u}$ and ${\bf v}$ that are neighbors on
${\mathbb C}^2(B)$, where ${\rm deg}({\bf u})$ is the number of
neighbors of ${\bf u}$, otherwise $p({\bf u,v}):=0$. 

Clearly when  $B=\{0\}$ we get the two-dimensional comb which inspired our 
choice for  the name of these particular anisotropic walks. We are interested 
in  the case when $B_n:=B\cap [ -n, n] , $ and  $|B_n| \sim c n^{\beta}$
with  some $0\leq\beta\leq1.$ Here, and in the sequel, $|B_n|$ stands for the 
(finite) number of elements in the set $B_n$. We are to discuss these comb-type 
walks at first in general, and then spell out three important cases, namely 
the case $\beta=1,$ $0<\beta<1,$ and $\beta=0.$  In each of these cases we 
prove strong approximation results for both components of the walks, by 
(time changed) Wiener processes. Our primary motivation was to connect and 
generalize some existing results on anisotropic walks. In each of these cases 
we use the fact that the number of the horizontal steps can be approximated via 
the occupation time of the set $B$ by the simple symmetric random walk on the 
$y-$axis defined by our approximation. It is important to emphasize that the 
scaling of the vertical component is always $N^{1/2},$  but the scaling of the 
horizontal component is $N^{1/2}$  only in  the case of $\beta=1.$ When $0<\beta<1,$ 
then the scaling depends on $\beta.$ Finally in the case of $\beta=0,$ the 
horizontal component is approximated by a time changed (by the local time of 
zero) Brownian motion resulting in scaling  by $N^{1/4.}$ 

In what follows we give some historical introduction, and try to  place our 
present paper in the context of some of the existing body of work on anisotropic walks. 
Initial studies of anisotropic walks are due to Silver {\it et al.} \cite{SSL},
Seshadri {\it et al.} \cite{SLS},  Shuler \cite{SH}, and Wescott \cite {WE}, 
who were motivated by the so called transport phenomena of statistical physics. 
Some of the most important contributions to the general anisotropic walk as in 
(\ref{prob1}) are due to Heyde \cite{H} and \cite{H93}, and Heyde {\it et al.} 
\cite{HWW}. 

As in  Heyde {\it et al.} \cite{HWW}, let  \begg
k^{-1}\sum_{j=1}^k p_{j}^{-1}=2\gamma_1+\ep_k,\quad \quad
k^{-1} \sum_{j=-k}^{-1} p_j^{-1}=2\gamma_2+\ep^*_k .   
\label{H1} \endd

\noindent
{\bf Theorem A} (\cite{HWW}) {\it For the anisotropic random  walk under 
condition {\rm(\ref{prob1})}, suppose that in {\rm (\ref{H1})} 
$\ep_k$ and $\ep^*_k$ are $o(1)$ as $k\to\infty.$ Then 
$$
\sup_{0\leq t\leq T} |n^{-1/2}C_2({[nt]})-Y(t)| \to 0  \quad a.s.
$$ 
as $n \to \infty,$ for all $T>0$,  where $\{Y(t), t \geq 0\} $ is a 
diffusion process on the same probability  space  as $\{C_2(n)\} $ whose 
distribution is defined by  
 $$Y(t)=W(A^{-1}(t)), \quad  t\geq 0,$$
where $\{W(t), t\geq 0 \}$ is a standard Brownian motion, $($or standard Wiener process$) $ and
 $$A(t)=\int_0^t\sigma^{-2}(W(s))\,ds$$
and}
 \begin{displaymath}
\sigma^2(y)=\left\{ 
\begin{array}{ll}
& \frac{1}{\gamma_1}\, \quad {\rm for}\quad y\geq 0, \\
\\
& \frac{1}{\gamma_2}\, \quad {\rm for}\quad y<0.
\end{array}\right.
\end{displaymath}

Here $A^{-1}(\cdot)$ is the inverse of $A(\cdot)$. The process $Y(t)$ is 
called oscillating Brownian motion if $\gamma_1\neq \gamma_2$, that is a 
diffusion with speed measure $m(dy)= 2\sigma^{-2}(y)dy.$

\medskip\noindent
{\bf Remark 1.1} Observe that $A(t)$ in the above theorem is equal to
\begg 
A(t)=\gamma_1\int_0^t I(W(s)\geq 0)\,ds+\gamma_2\int_0^t I(W(s)< 0)\,ds.
\label{at}
\endd 

In  an earlier paper Heyde {\cite H}
used the following somewhat more restrictive asymptotic  version of (\ref{H1}), 
when  
 \ $\gamma_1=\gamma_2=\gamma$: 
\begg n^{-1}\sum_{j=1}^n p_j^{-1}=2\gamma+o(n^{-\tau}),\qquad\qquad
n^{-1}\sum_{j=1}^n p_{-j}^{-1}=2\gamma+o(n^{-\tau}), \qquad
\label{bum}
\endd
as $n\to\infty$, for some constants $\gamma$,  $1<\gamma< \infty$ and
$1/2<\tau<\infty.$ 
Under (\ref{bum}) he proved a strong approximation result for the second 
coordinate by a  rescaled  Brownian motion.
Under the same condition the following  simultaneous strong 
approximation result was proved for $(C_1(\cdot), C_2(\cdot)).$

\medskip\noindent
{\bf Theorem B} (\cite{CCFR13}) {\it Under conditions {\rm(\ref{prob1})} and 
{\rm (\ref{bum})} with 
$1/2<\tau\leq 1$, on an appropriate probability space for the random walk
$$\{{\bf C}(N)=(C_1(N),C_2(N));\, \, N=0,1,2,\ldots\},$$
one can construct two independent
standard Wiener processes $\{W_1(t);\, t\geq 0\}$, $\{W_2(t);\,
t\geq 0\}$ so that, as $N\to\infty$, we have with any}
$\varepsilon>0$
\begg
\left|C_1(N)-W_1\left(\frac{\gamma-1}{\gamma}\,N\right)\right|+
\left|C_2(N)-W_2\left(\frac{1}{\gamma}\,N\right)\right|
=O(N^{5/8-\tau/4+\varepsilon})\quad a.s.
\endd

We note that Theorems A and B are true for $p_j$, more general than 
given in (\ref{combl}).

In \cite{CCFR13} and  \cite {CFR14} related  issues of  recurrence, local time and range  are 
discussed. Theorem 1.1  in \cite {CFR14}
 implies that all comb-type walks are recurrent.

Now we  mention some  particular cases of the random walk ${\bf C}(N)$ as in  
(\ref{prob1}) .

The case $p_j=1/4,\quad j=0, \pm1,\pm2....$ is the simple symmetric walk on 
the plane for which we refer to Erd\H{o}s and Taylor \cite{ET}, Dvoretzky   
and Erd\H{o}s \cite{DE} and R\'ev\'esz \cite {RE}. 
 
As we mentioned earlier, the case when  $p_0=1/4$ and $p_j=1/2$ with 
$j=\pm1,\pm2....$ means that all horizontal lines except the $x$-axis are 
missing. This is the so called random walk on the two dimensional comb, for 
which  $\gamma_1=\gamma_2=1, $ and hence is excluded from  Theorem B above. 
This model and some similar ones have many applications in physics, so a 
number of these early results are in the physics literature (\cite{A}, 
\cite{AR2}, \cite{AR4}, \cite{CAR}, \cite{DJ}, \cite{DJW}, \cite{REY}, 
\cite{ZA}, \cite{ZAE}). A present study of these models is provided in a new 
book by Iomin {\it et al.} \cite {IMH}. From the many papers about the exact 
comb  model as above, we only wish to mention the following few. Weiss and 
Havlin  \cite{WH} derived an asymptotic formula  for the $n$-step transition 
probability ${\bf P}({\bf C}(N)=(k,j))$, Bertacchi \cite {BE} was the first 
who noted that while a Brownian motion is the right object to approximate 
$C_2(\cdot)$, for the first component $C_1(\cdot),$ the right  approximation 
is by a Brownian motion time changed by the local time of the second component,
and she proved a simultaneous weak convergence result for the two components. 
Bertacchi and Zucca \cite{BZ} investigated the ratio of the $n$-step  vertical 
and horizontal transition probabilities, suggesting that this walk spends most 
of its time on some  teeth of the comb. In Cs\'aki {\it et al.} \cite{CCFR08} 
we established a simultaneous strong approximation for the two coordinates of 
the random walk ${\bf C}(N)=(C_1(N),C_2(N))$ that reads as follows.

\medskip\noindent
{\bf Theorem C} (\cite{CCFR08}) {\it On an appropriate probability space 
for the simple random walk 
\newline $\{{\bf C}(N)=(C_1(N),C_2(N)); N=0,1,2,\ldots\}$ on the 
two-dimensional comb lattice ${\mathbb C}^2,$ one can construct two 
independent standard Wiener processes $\{W_1(t);\, t\geq 0\}$, 
$\{W_2(t);\, t\geq 0\}$ so that, as $N\to\infty$, we have with any 
$\varepsilon>0$
$$
N^{-1/4}|C_1(N)-W_1(\eta_2(0,N))|+N^{-1/2}|C_2(N)-W_2(N)|
=O(N^{-1/8+\varepsilon})\quad a.s.,
$$
where $\eta_2(0,\cdot)$ is the local time process at zero of}
$W_2(\cdot)$.

Some Strassen type theorems are also given in the same paper. Furthermore, 
the local time of the comb walk is investigated in Cs\'aki {\it et al.} 
\cite{CCFR011} . 

In Cs\'aki {\it et al.} \cite{CCFR12} we investigated another special case, 
when $p_j=1/4,\, j=0,1,2,\ldots$, $p_j=1/2,\, j=-1,-2,\ldots, $ namely all 
horizontal lines under the $x$-axis are deleted, hence it is a simple random 
walk on the half-plane half-comb (HPHC) structure. In this case  $\gamma_1=2$ 
and $\gamma_2=1,$ thus Theorem B is not applicable. However this model 
satisfies the conditions of Theorem A and hence the second coordinate can be 
approximated by an oscillating Brownian motion. Our main result therein reads 
as follows.

\medskip\noindent
{\bf Theorem D} (\cite{CCFR12}) {\it On an appropriate probability space for 
the HPHC random walk \newline $\{{\bf C}(N)=(C_1(N),C_2(N)); N=0,1,2,\ldots\}$ 
with $p_j=1/4,\, j=0,1,2,\ldots$, $p_j=1/2,\, j=-1,-2,\ldots$ one can 
construct two independent standard Wiener processes 
$\{W_1(t);\, t\geq 0\}$, $\{W_2(t);\, t\geq 0\}$ such that, as $N\to\infty$, 
we have with any} $\varepsilon>0$
$$
|C_1(N)-W_1(N-A_2^{-1}(N))|+|C_2(N)-W_2((A_2^{-1}(N))|
=O(N^{3/8+\varepsilon})\quad a.s.,
$$
{\it where} 
$A_2(t)=2\int_0^t I(W_2(s)\geq 0)\,ds+\int_0^t I(W_2(s)< 0)\,ds$ \,
{\it as in} (\ref{at}).

In what follows, we are to 
generalize Theorem C and D under a somewhat relaxed version  of condition   
(\ref{bum}), when the $\gamma$ values in the two sums are different, discuss 
its consequences and consider some other interesting choices of the set $B.$ 

The structure of this paper from now on is as follows. In Section 2
we give preliminary facts and results. In  Section 3, first we redefine the 
walk on ${\mathbb C}^2(B)$ in terms of two independent simple symmetric walks. 
Then we list some facts which do not depend on the choice of $B,$ and prove 
some results which we will need for the rest of the paper. Section 4 contains 
our main results. In Section 5 some further questions and problems are
discussed.  

\section{Preliminaries}
\renewcommand{\thesection}{\arabic{section}} \setcounter{equation}{0}

\setcounter{thm}{0} \setcounter{lemma}{0}
In this section we list some well-known results, and some new ones which 
will be used in the rest of the paper. In case of the known ones  we won't 
give the most general form of the results, just as much as we intend to use,
while the exact reference will also be provided for the interested
reader.

Let $\{X_i\}_{i\geq1}$ be a sequence of independent i.i.d. random variables, 
with ${\mathbf P}(X_i=\pm1)=1/2.$ Then the simple symmetric random walk on 
the line is defined as $S(n)=\sum_{i=1}^n X_i, $  and its local time is  
$\xi(j,n)=\#\{k: 0< k \leq n,  S(k)=j\}, \,\,n=1,2...,$  for any integer $j.$ 

Define $M_n=\max_{0\leq k \leq n}|S(k)|.$  Then we have the usual law of the
iterated logarithm (LIL) and Chung's LIL \cite{CH}.

\medskip\noindent
{\bf Lemma A}
$$\limsup_{\noo} \frac{M_n}{(2n \log\log n)^{1/2}}=1, \qquad
\liminf_{\noo}\left(\frac{\log\log
n}{n}\right)^{1/2}M_n=\frac{\pi}{\sqrt{8}}\qquad a.s. $$
 
For $\xi(n)=\sup_{x}\xi(x,n)$ we have Kesten's LIL for local time.

\medskip\noindent {\bf Lemma B} (Kesten \cite{K}) {\it For the maximal local 
time we have}
$$ \limsup_{n\to\infty} \frac{\xi(n)}{(2n\log \log n)^{1/2}}=
1\quad a.s.$$

According to the lower lower class (LLC) result for the local time 
(see e.g. R\'ev\'esz \cite{RE}, page 119), the following holds true.

\medskip\noindent
{\bf Lemma C} (\cite{RE}) {\it  For the  local time of the simple
symmetric walk  we have for any $\varepsilon>0$} {\it and large enough} $n$
$$\xi(0,n)\geq \frac{\sqrt{n}}{(\log n)^{1+\ep}} \qquad a. s.$$

In Cs\'aki and F\"oldes \cite{CF87} the following stability result was
concluded for local time.

\medskip\noindent {\bf Lemma D} (\cite{CF87}) {\it  For the  local time of 
the simple symmetric random walk we have  that with any $\varepsilon>0$ and  }
  $\displaystyle{h(n)=
\frac{\sqrt{n}}{(\log n)^{1+\ep}}},$ 
$$\lim_{\noo}\sup_{|x|\leq h(n)}\left|\frac{\xi(x,n)}{\xi(0,n)}-1
\right|=0 \quad a.s.$$

In Cs\'aki and R\'ev\'esz \cite{CR83} the following result was given about the 
uniformity of the local time.

\medskip\noindent
{\bf Lemma E} (\cite{CR83}) {\it For the simple symmetric random walk for any 
$\varepsilon>0$ we have}
$$\lim_{n\to \infty}
\frac{\sup_{x\in \mathbb{Z}}|\xi(x+1,n)-\xi(x,n)|}{n^{1/4+\varepsilon}}=0\quad
a.s.
$$

\noindent
{\bf Remark 2.1} In fact, Lemma E deals with more general random walks, but 
we only need it for a simple symmetric random walk.

\medskip
The following result is a version of Hoeffding's inequality, which
is explicitly stated in T\'oth \cite{TO}.

\medskip\noindent{\bf Lemma F} (\cite{TO}) {\it Let $G_i$, $i=1,2,\ldots$ be 
i.i.d.random variables with the common geometric distribution 
${\bf P}(G_i=k)=2^{-k-1},\quad k=0,1,2...$ Then

$${\bf P}\left(\max_{1\leq j\leq n}
\left|\sum_{i=1}^j (G_i-1)\right|>\lambda\right)\leq
2\exp(-\lambda^2/{8n})
$$
for $0<\lambda<na$ with some} $a>0.$

Let $(W(t),\, t\geq 0)$ be a standard Wiener process .
 Its local time $(\eta(x,t),\, x\in \mathbb R,\, \, t\geq 0)$
 is defined as
$$
\eta(x,t)=\lim_{\varepsilon\to 0}\frac{1}{2\varepsilon}
\int_0^t I\{W(s)\in (x-\varepsilon,x+\varepsilon)\}\, ds,
$$
where $I\{\cdot\}$ denotes the indicator function.

Concerning the increments of the Wiener process we quote the following
result from Cs\"org\H o and R\'ev\'esz \cite{CSR79}, page 69.

\medskip\noindent
{\bf Lemma G} (\cite{CSR79}) {\it Let $0<a_T\leq T$ be a non-decreasing 
function of $T$. Then, as $T\to\infty$, we have}
$$
\sup_{0\leq t\leq T-a_T}\sup_{s\leq a_T}|W(t+s)-W(t)|=
O(a_T(\log(T/a_T)+\log\log T))^{1/2} \qquad a.s.
$$
The above statement is also true if $W(\cdot)$ replaced by the simple 
symmetric random walk $S(\cdot).$

We quote the following simultaneous strong approximation result from 
R\'ev\'esz \cite{RE81}. 

\medskip\noindent
{\bf Lemma H} (\cite{RE81}) {\it On an appropriate probability space for a 
simple symmetric random walk 

\noindent
$\{S(n);\, n=0,1,2,\ldots\}$ with local time
$\{\xi(x,n);\, x=0,\pm1,\pm2,\ldots;\, n=0,1,2,\ldots\}$ one can
construct a standard Wiener process $\{W(t);\, t\geq 0\}$ with local time
process $\{\eta(x,t);\, x\in\mathbb R; \, t\geq 0\}$ such that, as
$n\to\infty$, we have for any $\varepsilon>0$
$$|S(n)-W(n)|=O(n^{1/4+\varepsilon})\quad {a.s.}
$$
and
$$
\sup_{x\in\mathbb Z}|\xi(x,n)-\eta(x,n)|=O(n^{1/4+\varepsilon})
\quad {a.s.},
$$
simultaneously.}

Finally we recall the following lemma from Klass \cite{KL}.

\medskip\noindent
{\bf Lemma I}   \cite{KL} Let $\{E_n\}_{n\geq1}$
be an arbitrary sequence of events such that ${\mathbf P}(E_n \,i.o.)=1.$
Let $\{F_n\}_{n\geq1}$  be another arbitrary sequence of events that is 
independent  of $\{E_n\}_{n\geq1},$ and assume that for some $n_0>0,$  
${\mathbf P}(F_n)\geq c>0,$ for $n>n_0.$  Then we have 
${\mathbf P}(E_n F_n\, i.o.)> 0.$
  
\section{The general case}
\renewcommand{\thesection}{\arabic{section}} \setcounter{equation}{0}
\setcounter{thm}{0} \setcounter{lemma}{0}

First we are to redefine our random walk $\{{\mathbf C}(N);\,
N=0,1,2,\ldots\}$. It will be seen that the process described right below
is equivalent to that given in the Introduction.

To begin with, on a suitable probability space consider two independent
simple symmetric (one-dimensional) random walks $S_1(\cdot)$, and
$S_2(\cdot)$. We may assume that on the same probability space we have a
sequence of i.i.d. geometric random variables 
$\{G_i,\,i\geq 1\}$ which are independent from $S_1(\cdot)$, and
$S_2(\cdot),$ with
$$
\mathbf{P}(G_i=k)=\left(\frac{1}{2}\right)^{k+1},\,\, k=0,1,2,\ldots 
\label{geoseq}
$$

We now construct our walk $\mathbf{C}(N)$ as follows. We will take all
the horizontal steps consecutively from $S_1(\cdot)$ and all the vertical
steps consecutively from $S_2(\cdot).$ 
First we will take some horizontal steps from $S_1(\cdot)$, then as many 
vertical steps from $S_2(\cdot),$ as needed to get to a level belonging to 
$B,$ then again some horizontal steps from $S_1(\cdot)$ and so on. Now we 
explain how to get the number of horizontal steps on each occasion.
Consider our walk starting from the origin. If the origin is in $B$ then  
take $G_1$ horizontal steps from $S_1(\cdot).$ (Note that $G_1=0$ is possible 
with probability $1/2$). If the origin does not belong to $B$ then the walk 
moves vertically taking its steps from $S_2(\cdot)$ until it hits a level 
belonging to $B.$  If this happens at the level $j$ then it moves horizontally 
on level $j$ taking $G_1$ steps  from $S_1(\cdot)$. Then it takes  some 
vertical steps from $S_2(\cdot)$, to arrive again to a level belonging to $B,$
where the next $G_2$ horizontal steps from $S_1(\cdot),$ should be taken,
and so on. In general, whenever the walk arrives at the level $j\in B$  
then take some horizontal steps, the number of which is given by the next in 
line (first unused) geometric random variables $G_i$. 

Let now $H_N,\, V_N$ be the number of horizontal and vertical steps, 
respectively, from the first $N$ steps of the just described process. 
Consequently, $H_N+V_N=N$, and
$$
\left\{{\bf C}(N);\, N=0,1,2,\ldots\right\}=
\left\{(C_1(N),C_2(N));\, N=0,1,2,\ldots\right\}
$$
\begg
\stackrel{{d}}{=}\left\{(S_1(H_N),S_2(V_N));\, N=0,1,2,\ldots\right\},
\label{equ}
\endd
where $\stackrel{{d}}{=}$ stands for equality in distribution.

Now we introduce a few more notations. Let $\xi_2(\cdot,\cdot)$ denote the 
local time of $S_2(\cdot)$.
\begin{eqnarray}
M_1(N)&=&M_1(N,B)=\max_{0\leq k\leq N}|C_1(N)|,\\
M_2(N)&=&M_2(N,B)=\max_{0\leq k\leq N}|C_2(N)|,\\
H_N&=&H_N(B)=\#\{ k: \quad 1< k\leq N, \,\,\, C_1(k)\neq C_1(k-1)\},\\
V_N&=&V_N(B)= \#\{ k: \quad 1< k\leq N, \,\,\, C_2(k)\neq C_2(k-1)\},\\
D_2(V_N)&=&D_2(V_N,B)=\sum_{j\in B}\xi_2(j,V_N). \label{dset}
\end{eqnarray}

\noindent  Clearly  $M_1(N), M_2(N)$ are the absolute
maxima of the first and second coordinates.  $H_N$
and $V_N$ are the number of horizontal and vertical steps, respectively, in 
the first $N$ steps of ${\bf C}(\cdot).$  $D_2(V_N)$ is the 
occupation time of $B$ by $S_2(\cdot)$ in the first $N$ steps of 
${\bf C}(\cdot).$

By Lemma B we have for any $\ep>0$, any $y\in\mathbb{Z}$ and $N$ large 
enough: 

\begin{itemize}
\item ${\rm {\bf Fact\,1.}\quad} \xi_2(y,V_N)
\leq (1+\ep)(2V_N\log\log V_N)^{1/2}\leq(1+\ep) (2N\log\log N)^{1/2}
\quad a.s. $

Combining Lemmas C and D we easily get
\item $ {\rm {\bf Fact\,2.}\quad} \displaystyle{ \min_{|y|\leq 
\frac{V_N^{1/2}}{(\log
V_N)^{1+\ep}}} \xi_2(y,V_N)\geq  (1-\ep) \frac{V_N^{1/2}}{(\log V_N)^{1+\ep}}}
\quad a.s.$

Lemma A gives the following two facts:
\item $  {\rm {\bf Fact\,3.}\quad} M_2(N) 
\leq (1+\ep)(2V_N\log\log V_N)^{1/2}\leq(1+\ep) (2N\log\log
N)^{1/2} \quad a.s.$
\item $  {\rm {\bf Fact\,4.}\quad} 
\displaystyle{ M_2(N)\geq(1-\ep)\pi\left(\frac{V_N}{8 \log \log
V_N}\right)^{1/2}} \quad a.s.$
\end{itemize}

 Lemmas A and D also imply 
\begin{itemize}
\item $  {\rm {\bf Fact\,5.}\quad} M_1(N)\leq (1+\ep) (2 H_N\log \log
H_N)^{1/2} \quad a.s.$
\item $ {\rm {\bf Fact\, 6.} \quad}  M_1(N)
\geq (1-\ep) \pi\left( \frac{H_N}{8\log \log H_N}
\right)^{1/2} \quad a.s.$

\noindent
\item ${\rm {\bf Fact\, 7.} }\quad
H_N\geq \sum_{\left\{y\in B, \,|y|\leq
\frac{V_N^{1/2}}{(\log V_N)^{1+\ep}}\right\}}
 \xi_2(y,V_N)
 \quad a.s.$
\end{itemize}

Consider $A(t)$ defined by (\ref{at}) and let $\alpha(t)=A(t)-t$.
\begin{lemma}
Let $\gamma_1>\gamma_2\geq1.$ Then 
\begin{itemize}
\item
$A(t)  \quad {\rm and}\quad \alpha(t) $ {\rm are nondecreasing}
\item
$\gamma_2 t \leq A(t) \leq \gamma_1t   \quad {\rm and} \quad  
\frac{t}{\gamma_1} \leq A^{-1}(t)\leq \frac{t}{\gamma_2} . \label{hard} $
\end{itemize}
\end{lemma}
\noindent{\bf Proof.}  Observe that

\begg A(t)=\gamma_2 t+(\gamma_1-\gamma_2)\int_0^t I(W(s)\geq 0)\,ds  
\label{ate}   , \endd

\begg\alpha(t)=A(t)-t =(\gamma_2-1)t +
(\gamma_1-\gamma_2)\int_0^t I(W(s)\geq 0)\,ds.   \label{alfa}\endd
Then our first statement is obvious.   

From (\ref{ate}) we have that $$\gamma_2 t \leq A(t) \leq \gamma_1t   $$
and 
$$ A^{-1}(\gamma_2 t)\leq A^{-1}(A(t))=t\leq A^{-1}(\gamma_1 t),  $$
that, in turn, implies

\begg \frac{t}{\gamma_1}\leq A^{-1}(t)\leq \frac{t}{\gamma_2}. 
\label{hard2}  \endd $\Box$

For a simple random walk with local time $\xi(\cdot,\cdot),$ let 
\begg
\widehat{A}(n)=\gamma_1 \sum_{j=0}^{\infty} \xi(j,n)+\gamma_2 
\sum_{j=1}^{\infty} \xi(-j,n).
\label{an}
\endd

\begin{lemma} On a probability space as in Lemma H
$$|\widehat{A}(n)-A(n)|=O(n^{3/4+\ep})\quad\quad a.s.$$
\end{lemma}

\noindent{\bf Proof.}
According  to Lemma 5.3 of Bass and Griffin \cite {BG} we have 
$$\sup_{t\leq T}\sup _{|k|\leq t^{1/2+\ep}} \sup_{|k|\leq z\leq|k|+1}  
|\eta(z,t)-\eta(k,t)|=O(T^{1/4+\ep})  \qquad  a.s.,$$
where $\eta(x,t)$ is the local time of a Wiener process.
Using this fact, Lemmas A and H, we have for the difference
\begin{eqnarray}&&|\widehat{A}(n)-A(n)| 
\leq\gamma_1\left(\sum_{j=0}^{n}|\xi(j,n)-\eta(j,n)|
+\sum_{j=0}^{n}|\int_j^{j+1}\eta(x,n)\,dx-\eta(j,n)|\right)\\  \nonumber
&+&\gamma_2\left(\sum_{j=1}^{n}|\xi(-j,n)-\eta(-j,n)|
+\sum_{j=1}^{n}|\int_{-j}^{-j+1}\eta(x,n)\,dx-\eta(-j,n)|\right)  
=O(n^{3/4+\ep})  \qquad  a.s.\label{gond}
\end{eqnarray}
$\Box$

Let $D_2(V_N)$ be defined by (\ref{dset}).
\begin{lemma} For any $\ep>0$, as $N\to\infty$,
$$\max_{1\leq i\leq N} |H_i-D_2(V_i) |=O (|D_2(V_N)|^{1/2+\ep})\quad a.s.$$
\end{lemma}

\noindent
{\bf Proof.}  Recall that $H_N$ is the number of horizontal steps in our 
random walk ${\bf C}(N),\, N=0,1,2,\ldots$ As it was explained in the 
construction, horizontal steps only occur on levels belonging to $B.$ When 
the vertical walk arrives to such a level, it takes some horizontal steps, 
the number of which follows geometric distribution with expected value 1. 
Thus the total number of horizontal steps is the sum of $D_2(V_N)$ i.i.d. 
geometric random variables, with expected value 1. However this statement is 
slightly incorrect, as if the $N$-th step is a horizontal one, the 
corresponding last geometric random variable might remain truncated. Denote by 
$H_N^+$ the number of horizontal steps which includes all the steps of this 
last geometric random variable. Then 
  $$H_N^+ -D_2(V_N)=\sum_{j=1}^{D_2(V_N)}(G_j-1), $$
where $G_i, i=1,2...$ are i.i.d. geometric random variables as
in Lemma F. According to this lemma

$${\bf P}\left(\max_{1\leq j\leq k}
\left|\sum_{i=1}^j (G_i-1)\right|>\lambda\right)\leq
2\exp(-\lambda^2/{8k}).
$$
Selecting $\lambda=k^{1/2+\ep}$ we get by the Borel-Cantelli lemma that for all 
large $k$
\begg \max_{1\leq j\leq k}
\left|\sum_{i=1}^j (G_i-1)\right|\leq k^{1/2+\ep}\quad a.s.  
\label{trivi} \endd
We have to apply this for $k=|D_2(V_N)|,$  and conclude that
$$\max_{1\leq i\leq N} |H_i^+-D_2(V_i)|=O(|D_2(V_N)|^{1/2+\ep})\quad a.s.,$$
as $N\to\infty$. Now to finish the proof, it is enough to observe that
$$H_N^+ -H_N\leq \max_{i\leq N}G_j, $$
where $G_j$ are geometric random variables with parameter $1/2$. Thus for any 
$\varepsilon>0$
$$
P\left(\max_{j\leq N}G_j>N^{\varepsilon}\right)\leq 
N\left(\frac{1}{2}\right)^{N^{\varepsilon}}
$$
and hence by the Borel-Cantelli lemma
$$
H_N^+-H_N\leq N^{\varepsilon}\quad a.s.
$$
for all large $N$. 
 $\Box$

\section{Main results}
\renewcommand{\thesection}{\arabic{section}} \setcounter{equation}{0}

\setcounter{thm}{0} \setcounter{lemma}{0}

Assume that we have the construction described in Section 3 with 
independent simple symmetric random walks $S_1(\cdot)$, $S_2(\cdot)$
together with two  independent standard Wiener processes $W_1(\cdot)$, 
$W_2(\cdot)$ satisfying Lemma H, i.e. for $j=1,2$
$$
|S_j(n)-W_j(n)|=O(N^{1/4+\varepsilon}) \quad a.s.
$$
and
$$
\sup_{x\in\mathbb Z}|\xi_j(x,n)-\eta_j(x,n)|=O(n^{1/4+\varepsilon})
\quad a.s.
$$

Now we define $A_2(\cdot)$ as in (\ref{at}), with $W$ replaced by $W_2$,
$\alpha_2(t)=A_2(t)-t$,
and also $\widehat{A}_2(\cdot)$ as in (\ref{an}), with $\xi$ replaced by 
$\xi_2$.

In what follows we assume that for $ B_n:=B\cap [-n, n]$ we have
\begg  |B_n|\sim c n^{\beta} \label{beta} \endd
with  some $0\leq\beta\leq1$ and with some constant $c>0.$ 
As we mentioned in the Introduction, 
we will separately consider the the cases $\beta=1$, $0<\beta<1$ and 
$\beta=0.$ 

\subsection{The case $\beta=1$}

To consider the case $\beta=1$, we now suppose that 
\begg n^{-1}\sum_{j=1}^n p_j^{-1}=2\gamma_1+o(n^{-\tau}), \qquad
n^{-1}\sum_{j=1}^n p_{-j}^{-1}=2\gamma_2+o(n^{-\tau}) \label{h}  \endd
as $n\to\infty$ with some $1/2<\tau\leq 1.$ 
Clearly $1\leq \gamma_1\leq 2$  and  $1\leq \gamma_2\leq 2.$ 

\begin{thm}Under the conditions {\rm (\ref{combl}), (\ref{h})} and 
$\max(\gamma_1,\gamma_2)>1$, on an appropriate probability space for the 
random walk $\{{\bf C}(N)=(C_1(N),C_2(N));\, \, N=0,1,2,\ldots\}$ one can 
construct two independent standard Wiener processes $\{W_1(t);\, t\geq 0\}$, 
$\{W_2(t);\, t\geq 0\}$ so that, as $N\to\infty$, we have with any
$\varepsilon>0$
\begg
\left|C_1(N)-W_1\left(N-A_2^{-1}(N)\right)\right|+
\left|C_2(N)-W_2\left(A_2^{-1}(N)\right)\right|
=O(N^{5/8-\tau/4+\varepsilon})\quad a.s. \label{hat}
\endd
where $A_2^{-1}(\cdot)$ is the inverse of $A_2(\cdot).$
\end{thm}

\noindent
{\bf Proof.} In what follows we take some ideas from Heyde \cite{H}.   
By our Lemma 3.3 
\begin{eqnarray} 
H_N&=&D_2(V_N)+O(N^{1/2+\ep})
=\sum_{j \in B}\xi_2(j,V_N)+O(N^{1/2+\varepsilon}) \\ \nonumber
&=&\sum_{j }\xi_2(j,V_N)\frac{1-2p_j}{2p_j}+O(N^{1/2+\varepsilon})
 =-V_N+\frac{1}{2}\sum_{j } \xi_2(j,V_N)\frac{1}{p_j}
+O(N^{1/2+\varepsilon})     \qquad  a.s..\label{fontos}
\end{eqnarray}
In the second line above we changed the summation for all $j,$ recalling  
the fact that $p_j=1/2$ for $j\notin B$, while  $p_j=1/4$ for $j\in B.$

To proceed, introduce the notation
$$\frac{1}{j}\sum_{k=1}^j \frac{1}{p_k}=\kappa_j, \quad
\frac{1}{j}\sum_{k=1}^j \frac{1}{p_{-k}}=\nu_j.$$

\begg N=H_N+V_N=\frac{1}{2}\sum_j \xi_2(j,V_N)
\frac{1}{p_j}+O(N^{1/2+\varepsilon})   \qquad  a.s.. \label{azN} \endd
Then
$$\sum_j \xi_2(j,V_N)\frac{1}{p_j}=\sum_{j=1}^\infty
\xi_2(j,V_N)(j\kappa_j-(j-1)\kappa_{j-1})+\sum_{j=1}^\infty
\xi_2(-j,V_N)(j\nu_j-(j-1)\nu_{j-1})+\xi_2(0,V_N)\frac{1}{p_0}$$
$$=\sum_{j=1}^\infty j \kappa_j(\xi_2(j,V_N)- \xi_2(j+1,V_N))+
\sum_{j=1}^\infty j \nu_j(\xi_2(-j,V_N)- \xi_2(-j-1,V_N))+
\xi_2(0,V_N)\frac{1}{p_0}$$
$$=\sum_{j=1}^\infty j(\kappa_j-2 \gamma_1)(\xi_2(j,V_N)-
\xi_2(j+1,V_N))+2\gamma_1 \sum_{j=1}^\infty j(\xi_2(j,V_N)- \xi_2(j+1,V_N))$$
$$+\sum_{j=1}^\infty j(\nu_j-2\gamma_2)(\xi_2(-j,V_N)-
\xi_2(-j-1,V_N))+2\gamma_2 \sum_{j=1}^\infty j(\xi_2(-j,V_N)-
\xi_2(-j-1,V_N))+\xi_2(0,V_N)\frac{1}{p_0}$$
$$=2\gamma_1 \sum_{j=0}^{\infty} \xi_2(j,V_N)+2\gamma_2 
\sum_{j=1}^{\infty} \xi_2(-j,V_N)+
\sum_{j=1}^{\infty} j(\kappa_j-2 \gamma_1)(\xi_2(j,V_N)
-\xi_2(j+1,V_N))$$
$$+\sum_{j=1}^{\infty}
j(\nu_j-2\gamma_2)(\xi_2(-j,V_N)- \xi_2(-j-1,V_N))
+\xi_2(0,V_N)\left(\frac{1}{p_0}-2\gamma_1\right).$$

Observe that from (\ref{h}) we have that
$$|j(\kappa_j-2 \gamma_1)|\leq c j^{1-\tau},\quad 
\quad |j(\nu_j-2\gamma_2)|\leq c |j|^{1-\tau}$$
for some $c>0.$ Now applying Lemma A for $S_2(\cdot)$, and Lemma E, we get
that
$$\sum_{j=1}^{\infty} j(\kappa_j-2 \gamma_1)(\xi_2(j,V_N)- \xi_2(j+1,V_N))
+\sum_{j=1}^{\infty } j(\nu_j-2\gamma_2)(\xi_2(-j,V_N)- \xi_2(-j-1,V_N))$$
$$=O(N^{1/4+\epsilon})\sum_{j=1}^{\max_{k\leq N}|S_2(k)|} j^{1-\tau}=
O(N^{1/4+\epsilon})O(N^{1-\tau/2+\ep})=O(N^{5/4-\tau/2+\ep}) \qquad  a.s.,$$
where here and throughout the paper the value of $\varepsilon$ might change
from line to line.

So we conclude, using Lemma 3.2, that
$$\frac{1}{2}\sum_j \xi_2(j,V_N)\frac{1}{p_j}=\gamma_1 
\sum_{j=0}^{\infty} \xi_2(j,V_N)+\gamma_2 \sum_{j=1}^{\infty} 
\xi_2(-j,V_N)+O(N^{5/4-\tau/2+\ep})$$
$$=\widehat{A}_2(V_N)+O(N^{5/4-\tau/2+\ep})=A_2(V_N)+O(N^{5/4-\tau/2+\ep})
+O(N^{3/4+\ep})
$$
$$=A_2(V_N)+O(N^{5/4-\tau/2+\ep}) \qquad  a.s.$$
since $1/2<\tau\leq1.$
Consequently from (\ref{azN})
$$N=A_2(V_N)+O(N^{5/4-\tau/2+\ep}) \qquad  a.s.$$
and 
$$V_N=A_2^{-1}(N)+O(N^{5/4-\tau/2+\ep}) \qquad  a.s.$$
\noindent
{\bf Remark 4.1} In the previous line we used the fact that 
$A_2^{-1}(u+v)-A_2^{-1}(u)\leq v.$  
To see this, first recall from Lemma 3.1 that $A_2(t),\,$  $A_2^{-1} (t)$  and  
$\alpha(t) =A_2(t)-t$ are all nondecreasing. Then
$$v=A_2(A_2^{-1}(u+v))-A_2(A_2^{-1}(u))=\alpha(A_2^{-1}(u+v))+A_2^{-1}(u+v)
-\alpha(A_2^{-1}(u))-A_2^{-1}(u)$$
$$\geq A_2^{-1}(u+v)-A_2^{-1}(u),$$

So we can conclude, using Lemmas H and G that 
$$C_2(N)=S_2(V_N)=W_2(V_N)+O(N^{1/4+\ep})=W_2((A_2^{-1}(N)
+O(N^{5/4-\tau/2+\ep}))$$
$$=W_2((A_2^{-1}(N))+O(N^{5/8-\tau/4+\ep}) \qquad  a.s.$$
Furthermore, by Lemmas H and G again 
$$C_1(N)=S_1(H_N)=S_1(N-V_N)=W_1(N-V_N)+O(N^{1/4+\ep})
$$
$$
=W_1(N-A_2^{-1}(N))+O(N^{5/8-\tau/4+\ep})) \qquad  a.s.,$$
proving our theorem. $\Box$

\noindent
{\bf Remark 4.2} Here we would like to mention that Theorem 4.1 is a 
generalization of Theorem D and Corollary 4.1 below contains  Corollary 4.4 
of \cite{CCFR12}.

\begin{corollary} Suppose that $\gamma_1>\gamma_2\geq 1.$ The following laws 
of the iterated logarithm hold.
\begin{itemize}
\item{\rm (i)}
$$
\limsup_{t\to\infty}\frac{W_1(t-A_2^{-1}(t))}{(t\log\log t)^{1/2}}=
\limsup_{N\to\infty}\frac{C_1(N)}{(N\log\log N)^{1/2}}
=\sqrt{2\left(1-\frac{1}{\gamma_1}\right)}
\quad a.s.,$$
\item{\rm (ii)}
$$\liminf_{t\to\infty}\frac{W_1(t-A_2^{-1}(t))}{(t\log\log t)^{1/2}}=
\liminf_{N\to\infty}\frac{C_1(N)}{(N\log\log N)^{1/2}}
=-\sqrt{2\left(1-\frac{1}{\gamma_1}\right)}
\quad a.s.,$$
\item{\rm (iii)}
$$
\limsup_{t\to\infty}\frac{W_2(A_2^{-1}(t))}{(t\log\log t)^{1/2}}=
\limsup_{N\to\infty}\frac{C_2(N)}{(N\log\log N)^{1/2}}=\sqrt{\frac{2}{\gamma_1}}
\quad a.s.,
$$
\item{\rm (iv)}
$$
\liminf_{t\to \infty}\frac{W_2(A_2^{-1}(t))}{(t\log\log t)^{1/2}}=
\liminf_{N\to\infty}\frac{C_2(N)}{(N\log\log N)^{1/2}}=-\sqrt{\frac{2}{\gamma_2}}
\quad a.s.
$$
\end{itemize}
\end{corollary}

\bigskip\noindent
{\bf Proof of (i) and (ii)}. By the law of the iterated logarithm for 
$W_1,$ and (\ref{hard2}) we have for all large enough $t$
$$
W_1(t-A_2^{-1}(t))\leq (1+\varepsilon)
(2(t-A_2^{-1}(t))\log\log(t-A_2^{-1}(t)))^{1/2}
$$
$$
\leq
(1+\varepsilon)\left(2t\left(1-\frac{1}{\gamma_1}\right)\log\log t 
\right)^{1/2}  \qquad  a.s.,
$$
which gives an upper bound in (i).

To give a lower bound in (i), for any sufficiently small $\delta>0$ define 
the events
$$
A_n^*=\{W_1(u_n)\geq (1-\delta)(2u_n\log\log u_n)^{1/2}\},
\quad B_n^*=\left\{\int_0^{\frac{u_n(1+\delta)}{(\gamma_1-1)}} 
I(W_2(s)\geq 0)\,ds>\frac{u_n}{\gamma_1-1}\right\},
$$
$n=1,2,\ldots$ Then, with some sequence $\{u_n\}$ ($u_n=a^n$ with 
sufficiently large $a$ will do), we have
$$
{\mathbf P}(A_n^*\, \, i.o.)=1,\qquad {\mathbf P}(B_n^*)>c>0.
$$
It follows from Lemma I that
$$
{\mathbf P}(A_n^*B_n^*\, \, \, \, i.o.)\geq c>0.
$$
By the 0\,-1 law this probability is equal to 1.
Recall that $\alpha_2(t)=A_2(t)-t$. We claim that if $B_n^*$ occurs then 
\begg \alpha_2\left(\frac{u_n(1+\delta)}{\gamma_1-1}\right)\geq u_n. 
\label{vegre}\endd
Now if $B_n^*$ occurs then by (\ref{alfa})
$$\alpha_2\left(\frac{u_n(1+\delta)}{\gamma_1-1}\right)
=\frac{(\gamma_2-1)u_n(1+\delta)}{\gamma_1-1}+(\gamma_1-\gamma_2)
 \int_0^{\frac{u_n(1+\delta)}{\gamma_1-1}} I(W_2(s)\geq 0)\,ds$$
 $$
 >\frac{(\gamma_2-1)u_n(1+\delta)}{\gamma_1-1}
+\frac{(\gamma_1-\gamma_2)u_n}{\gamma_1-1}\geq u_n.$$
Let $t_n$ be defined by
$$
u_n=t_n-A_2^{-1}(t_n)=\alpha_2(A_2^{-1}(t_n)).
$$
Since $B_n^*$ implies
$$
\alpha_2\left(\frac{u_n(1+\delta)}{\gamma_1-1}\right)>u_n
=\alpha_2(A_2^{-1}(t_n)),
$$
and $\alpha_2(\cdot)$ is nondecreasing, we  have that
$$\frac{u_n(1+\delta)}{\gamma_1-1}>A_2^{-1}(t_n).$$
Thus, using (\ref{hard2})

$$u_n>\frac{\gamma_1-1}{1+\delta}A_2^{-1}(t_n)>\frac{\gamma_1-1}{1+\delta} 
\frac{t_n}{\gamma_1}=(1-\frac{1}{\gamma_1}) \frac{t_n}{1+\delta}.$$
Hence $A_n^*B_n^*$ implies
$$
W_1(t_n-A_2^{-1}(t_n))\geq (1-\delta)\left(
\frac{2(1-\frac{1}{\gamma_1}) t_n\log\log t_n}{1+\delta}\right)^{1/2}  
\qquad  a.s.
$$
Since $\delta>0$ is arbitrary, this gives a lower bound in (i).

The proof of (ii) follows by symmetry.

\noindent
{\bf Proof of (iii)}. 
We have infinitely often with probability 1
$$
W_2(A_2^{-1}(t))\geq (1-\varepsilon)(2A_2^{-1}(t)\log\log t)^{1/2}
\geq(1-\varepsilon)\sqrt{\frac{2}{\gamma_1}}(t\log\log t)^{1/2},
$$
where we used  (\ref{hard2})  to get the second inequality. 

To give an upper bound, we use the formula for the distribution of
the supremum of $W_2(A_2^{-1}(t))$ given in Corollary 2 of Keilson and 
Wellner \cite{KW}, which in our case is equivalent to
$$
{\mathbf P}(\sup_{0\leq s\leq t}W_2(A_2^{-1}(s))>y)
$$
$$
=\frac{4\sqrt{\gamma_1}}{\sqrt{\gamma_1}+\sqrt{\gamma_2}}
\sum_{k=0}^\infty
\left(\frac{\sqrt{\gamma_2}-\sqrt{\gamma_1}}{\sqrt{\gamma_1}
+\sqrt{\gamma_2}}\right)^k
\left(1-\Phi\left(\frac{(2k+1)\sqrt{\gamma_1}}{\sqrt{t}}y\right)\right) 
\qquad y\geq 0.
$$
From this it is easy to give the estimation
$$
{\mathbf P}\left(\sup_{0\leq s\leq t}W_2(A_2^{-1}(s))>y\right)\leq
c\exp\left(-\frac{\gamma_1y^2}{2t}\right)
$$
with some constant $c$, from which the upper estimation in (iii) follows 
by the usual procedure.
{\bf Proof of (iv)}. The lower estimation is easy. Namely, by (\ref{hard2}) 
we have for $t$ big enough that
$$
W_2(A_2^{-1}(t))\geq -(1+\varepsilon)(2A_2^{-1}(t)\log\log A_2^{-1}(t))^{1/2}
\geq -(1+\varepsilon)\left(2\frac{t}{\gamma_2}\log\log t\right)^{1/2} 
\quad  a.s.
$$
It remains to prove an upper estimation in (iv). By the law of the 
iterated logarithm for $W_2(\cdot)$ 
\begin{equation}
W_2(v)\leq -((2-\varepsilon)v\log\log v)^{1/2}
\label{lil}
\end{equation}
almost surely for infinitely many $v$ tending to infinity. Define
\begg \mu(v)=\int_0^v I(W_2(s)\geq 0)\,ds.  
\label {most} \endd
Let $\zeta(v)$ be the last zero of $W_2(\cdot)$ before $v$, i.e.,
$$
\zeta(v)=\max\{u\leq v:\, W_2(u)=0\}, 
$$
By Theorem 1 of Cs\'aki and Grill \cite{CG}, (or from the Strassen theorem), 
for large $v$ satisfying 
(\ref{lil}) we have $\zeta(v)\leq \varepsilon v$, and hence also 
$\mu(v)\leq\zeta(v)\leq \varepsilon v$. Now put $v=A_2^{-1}(t)$, i.e., 
$t=A_2(v)\leq \ep v\gamma_1+\gamma_2 v=(\ep \gamma_1+\gamma_2 )v$, from which 
$v=A_2^{-1}(t)\geq t/(\ep \gamma_1+\gamma_2 )$. Hence
$$
W_2(v)=W_2(A_2^{-1}(t))\leq 
-\left((2-\varepsilon)\frac{t}{\ep \gamma_1+\gamma_2 }\log\log t\right)^{1/2}
$$  
infinitely often with probability 1.
Since $\varepsilon>0$ is arbitrary, this gives an upper bound in (iv). $\Box$

It is not hard to calculate the density function of $A_2^{-1}(t)$ 
and $t-A_2^{-1}(t).$
\begin{lemma} Suppose that $\gamma_1>\gamma_2\geq 1.$ 
$${\mathbf P}(A_2^{-1}(t)\in dv) =\frac{t}{\pi v}\frac{1}
{\sqrt{(v\gamma_1-t)(t-\gamma_2v)}}\,dv
\quad {\rm for}\quad \frac{t}{\gamma_1}<v<\frac{t}{\gamma_2},$$
$${\mathbf P}(t-A_2^{-1}(t)\in dv) 
=\frac{t}{\pi(t-v)}\frac{1}{\sqrt{((\gamma_1-1)t-\gamma_1 v)(t(1-\gamma_2) 
+\gamma_2v)}}\,dv$$
$${\rm for}\quad t 
\left(1-\frac{1}{\gamma_2}\right)<v<t\left(1-\frac{1}{\gamma_1}\right).$$
\end{lemma}  
{\bf Proof.}  Recall the definition of $\mu(v)$ from (\ref{most}).
\begin{eqnarray}
{\mathbf P}(A_2^{-1}(t)<v)&=&{\mathbf P}(t<A_2(v))=
{\mathbf P}(t<\gamma_2 v+(\gamma_1-\gamma_2)\mu(v))\\ 
\nonumber
&=&P\left(\frac{\mu(v)}{v}>\frac{t-\gamma_2v}{v(\gamma_1-\gamma_2)}\right)= 
1-\frac{2}{\pi}\arcsin\left(\sqrt{\frac{t-v\gamma_2}{v(\gamma_1-\gamma_2)}}\,
\right). \nonumber
\end{eqnarray}
By differentiation we get the first statement  and the second goes similarly. 
$\Box$

\subsection{The case $0<\beta<1$}

Now we want to consider the second case when we have $0<\beta<1$ in 
(\ref{beta}), which means that a considerable portion of the horizontal lines 
are missing. Recall that 
$$D_2(V_N)=\sum_{j\in B}\xi_2(j,V_N). $$
\begin{thm} Under the condition {\rm (\ref{beta})} with $0<\beta<1,$
on an appropriate probability space for the random walk $\{{\bf
C}(N)=(C_1(N),C_2(N)); N=0,1,2,\ldots\}$  one can construct two independent
standard Wiener processes $\{W_1(t);\, t\geq 0\}$, $\{W_2(t);
t\geq 0\}$ so that, as $N\to\infty$, we have with any
$\varepsilon>0$
$$
\left|C_1(N)-W_1\left(D_2(V_N\right)\right)|=O(N^{1/8+\beta/8+\varepsilon})
\quad a.s.$$
$$|C_2(N)-W_2(N)|
=O(N^{1/4+\beta/4+\ep})   \quad a.s.$$
\end{thm}

\noindent
{\bf Proof.} Using Lemmas A and B and (\ref{beta}) we conclude
$$D_2(V_N)=O(N^{1/2 +\beta/2+\ep})\quad a.s.$$
Clearly from Lemma 3.3
$$H_N=D_2(V_N)+ O(D_2(V_N))^{1/2+\ep}=D_2(V_N)+O(N^{1/4+\beta/4+\ep})
\quad a.s.$$
and by Lemmas H and G
$$C_1(N)=S_1(H_N)=W_1(H_N)+O(H_N^{1/4+\ep})=W_1(H_N)
+O( N^{1/8+\beta/8+\ep})$$
$$=  W_1(D_2(V_N))+O( N^{1/8+\beta/8+\ep}) \quad a.s.$$
$$C_2(N)=S_2(V_N)=W_2(V_N)+O(N^{1/4+\ep})=W_2(N-H_N)
+O(N^{1/4+\ep})$$
$$=W_2(N)+O(N^{1/4+\beta/4+\ep})+O(N^{1/4+\ep})
=W_2(N)+O(N^{1/4+\beta/4+\ep}) \quad a.s.$$
$\Box$

\noindent
{\bf Remark 4.2} It is very easy to see  using Lemmas G and H, that the first 
statement in Theorem  4.2 can be replaced by the following more natural one: 

$$
\left|C_1(N)-W_1\left(L_2(V_N\right)\right)|=O(N^{1/8+\beta/4+\varepsilon})
\quad a.s.$$
where $L_2(V_N):=\sum_{j\in B}\eta_2(j,V_N),$  that is to say using the  
Wiener local time instead of the random walk local time. However this change 
slightly weakens the rate  of the approximation. Moreover we can replace  
$ L_2(V_N)$ in the argument  of $W_1(\cdot)$ with $L_2(N)$ using the local 
time increment result in Cs\'aki {\it et al.}  \cite{CCFR83} (see e.g. 
\cite{RE}, page 121, Theorem 11.9), to get the following  weaker approximation:

$$
\left|C_1(N)-W_1\left(L_2(N \right)\right)|=O(N^{1/8+3\beta/8+\varepsilon})
\quad a.s.$$

Clearly, the weakness of Theorem 4.2 is, that we do not know the limiting 
distribution and other limit theorems (lower and upper classes) for $D_2(V_N).$
 So we won't be able to get LIL-s for $C_1(N)$ as in the case of Theorem 4.1. 
In what follows we get some simple estimates instead.
Namely, using Lemmas A and B and  (\ref{beta}),  we have 
that with some positive constant $c>0$
$$D_2(V_N)\leq c (N\log\log N)^{\frac{1+\beta}{2}} \quad  a.s.$$

On the other hand, from Lemmas C and D and (\ref{beta})  we have for  any
$\ep>0$ that
$$D_2(V_N)\geq \frac{N^{\frac{1+\beta}{2}}}{(\log N)^{1+\beta+\ep}} \quad  
a.s.$$
Using these bounds we can conclude that with some $c_1>0$
$$\frac{N^{\frac{1+\beta}{4}}}{(\log N)^{\frac{1+\beta}{2}+\ep}}
\leq C_1(N)\leq c_1 N^{\frac{1+\beta}{4}}(\log \log N)^{\frac{3+\beta}{4}}
\quad  a.s.$$

Furthermore, for the second coordinates of our walk we have 
\begg
\limsup_{N\to\infty} \frac{\max_{0\leq k\leq N}|C_2(k)|}
{(2N\log\log N)^{1/2}}=1\qquad a.s. \endd
From  the so-called other law of the iterated
logarithm due to Chung \cite{CH} we obtain for $C_2(N)$ 
\begin{equation}
\liminf_{N\to\infty}\left(\frac{8\log\log N}{\pi^2 N}\right)^{1/2}
\max_{0\leq k\leq N}|C_2(k)|=1 \qquad a.s.
\end{equation}

\subsection{The case $\beta=0$}

Now we consider the case when $B$ is finite, that is to say, we only have 
finitely many horizontal lines (so in (\ref{beta}) $\beta=0).$ 
\begin{thm}  
Suppose that $B$ is finite. Then on an appropriate probability space for the 
random walk  $\{{\bf C}(N)=(C_1(N),C_2(N)); N=0,1,2,\ldots\}$ one can 
construct two independent standard Wiener processes 
$\{W_1(t);\, t\geq 0\}$, $\{W_2(t);
t\geq 0\}$ so that, as $N\to\infty$, we have with any $\varepsilon>0$
$$
N^{-1/4}|C_1(N)-W_1(|B|\eta_2(0,N))|+N^{-1/2}|C_2(N)-W_2(N)|
=O(N^{-1/8+\ep})\quad a.s.$$
\end{thm}

\noindent
{\bf Proof.} By Lemma 3.3 and Lemmas B, E and H

$$H_N=D_2(V_N)+(D_2(V_N))^{1/2+\ep}=|B|\xi_2(0,N)+O(N^{1/4+\ep})
 =|B|\eta_2(0,N)+O(N^{1/4+\ep})\quad a.s.$$
Hence by Lemmas H and  G we have
\begin{eqnarray}
C_1(N)&=&S_1(H_N)=S_1(|B|\eta_2(0,N)+O(N^{1/4+\ep})) \\ \nonumber
&=&W_1(|B|\eta_2(0,N)+O(N^{1/4+\ep}))+O(N^{1/8+\ep})
=W_1(|B|\eta_2(0,N) )+O(N^{1/8+\ep}) \quad a.s.\\ \nonumber
\end{eqnarray}
Similarly
$$C_2(N)=S_2(V_N)=S_2(N)+O(H_N^{1/2+\ep})=W_2(N)+O(N^{1/4+\ep})\quad a.s.$$ 
$\Box$

Observe that Theorem 4.3 is a generalization of Theorem C. Now we list some 
conclusions of Theorem 4.3.

Define the continuous version  of the random walk process on our lattice, 
having horizontal lines only in $B$ as follows:
$$\{{\bf C}(xN)=(C_1(xN),C_2(xN)): 0\leq x\leq1\}.$$  
We have almost surely, as $N\to\infty$,
$$ \sup_{0\leq x\leq 1}\left\Vert\left(\frac{C_1(xN)-W_1(|B|\eta_2(0,xN))}
{N^{1/4}(\log\log N)^{3/4}},
\frac{C_2(xN)-W_2(xN)}{(N\log \log N)^{1/2}}\right)\right\Vert\to 0,$$
where $||.||$ means the Euclidean norm.
We have the following laws of the iterated logarithm (for the first component 
see Theorem 2.2 in Cs\'aki {\it et al.}  \cite{CCFR95}).

$$ \limsup_{n\to \infty} \frac{C_1(N)}{\sqrt{|B|}N^{1/4}(\log\log N)^{3/4}}
=\frac{2^{5/4}}{3^{3/4}}  \quad a.s. \quad {\rm and} \quad
\limsup_{N\to \infty} \frac{C_2(N)}{(2N\log\log N)^{1/2}}=1 \quad a.s.$$

As to the liminf behavior of the max functionals of
the two components,  we have the same results as for the two dimensional comb 
lattice \cite{CCFR08}. These results are based on the corresponding ones  
for Wiener process  and  the iterated process $W_1(\eta_2(0,t))$ and  
the work of Chung \cite {CH}, Hirsch \cite{HI}, Bertoin \cite{BER}, and 
Nane \cite{NA} . 

Based on \cite{NA}, we get the following: Let $\rho(n),\, n=1,2,\ldots$, be 
a non-increasing sequence of positive numbers such that $n^{1/4}\rho(n)$ is
non-decreasing. Then we have almost surely that

$$
\liminf_{N\to\infty}\frac{\max_{0\leq k\leq
N}C_1(k)}{N^{1/4}\rho(N)}=0\quad or\quad \infty
$$
and

$$
\liminf_{n\to\infty}\frac{\max_{0\leq k\leq
N}C_2(k)}{N^{1/2}\rho(N)}=0\quad or\quad \infty,
$$
according as to whether the series $\sum_1^\infty \rho(n)/n$ diverges or
converges.

\begin{equation}
\liminf_{N\to\infty}\left(\frac{8\log\log N}{\pi^2 N}\right)^{1/2}
\max_{0\leq k\leq N}|C_2(k)|=1 \qquad a.s.
\end{equation}

On the other hand, for the max functional of $|C_1(\cdot)|$  we obtain 
from \cite{CCFR08} the following result.

Let $\rho(n),\, n=1,2,\ldots$, be a non-increasing
sequence of positive numbers such that $n^{1/4}\rho(n)$ is
non-decreasing. Then we have almost surely that

$$
\liminf_{N\to\infty}\frac{\max_{0\leq k\leq
N}|C_1(k)|}{N^{1/4}\rho(N)}=0\quad or\quad \infty,
$$
as to whether the series $\sum_1^\infty \rho^2(n)/n$ diverges 
or converges.

\section{ Final  remarks  and  further questions }     
\renewcommand{\thesection}{\arabic{section}} \setcounter{equation}{0}

\setcounter{thm}{0} \setcounter{lemma}{0}

As we mentioned in the Introduction, the main advantage of using  
(\ref{combl}) instead of (\ref{prob1}) in this paper, is that the number of horizontal steps 
can be easily approximated by the total time spent by the vertical walk in 
the set $B,$ as  shown in Lemma 3.3. It would be interesting to get the 
results of this paper under the much less restrictive  transition probabilities of 
(\ref{prob1}). Some further areas of investigation would be to discuss the 
other usual questions on random walks on the plane. Namely we would like to 
investigate the return time to zero,  local times and the range. These issues 
were considered in \cite {CCFR13} only for the so called periodic case, when 
$p_j=p_{j+L}$ for each $j\in \mathbb{Z},$ where $L\geq 1$ is a positive integer.

To illustrate our results we would like to mention some examples. 

\smallskip
\noindent
{\bf Example 1.} In  \cite{CCFR13}   we discussed a special periodic case,  
the so called uniform case $p_j=1/4$ if $|j|=0(mod\, L )$ and $p_j=1/2$ 
otherwise. Now our Theorem 4.1 contains the case where 

\begin{eqnarray*}
p_j&=&1/4 \quad {\rm if}\quad j=0(mod\, L )\quad {\rm for}\quad j\geq 0\\
p_j&=&1/4 \quad {\rm if}\quad j=0(mod \,K )\quad {\rm for}\quad j< 0\\
p_j&=&1/2  \quad  {\rm otherwise.}
\end{eqnarray*}
Then we get $\gamma_1=\frac{L+1}{L}$ and $\gamma_2=\frac{K+1}{K}$ for any 
pairs of integers $L$ and $K$ in (\ref{h}) and Theorem 4.1 holds with 
$\tau=1.$ Clearly, many more elaborate patterns of keeping and discarding  
horizontal lines could be handled.  

\smallskip
\noindent
{\bf Example 2.} As we mentioned earlier, the case $\gamma_1=\gamma_2$ in 
Theorem 4.1 gives our Theorem B. On the other hand, Theorem D  about the HPHC 
walk means that in Theorem 4.1 we have $\gamma_1=2,$ and $\gamma_2=1$.  But 
Theorem 4.1 permits many HPHC like structures as the  following example shows. 
For any $\alpha_1>1, \alpha_2>1$ let 
$$B^{(1)}=\{j: j=0,1,2,...\}, \quad B^{(2)}=\{[k^{\alpha_1}], k=1,2,...\},
\quad B^{(3)}=\{- [k^{\alpha_2}], k=1,2,..  \}$$
and 
$$B=(B^{(1)}\setminus B^{(2)})\cup B^{(3)},$$
namely above the $x-$axis, which is the plane part, we eliminate many horizontal 
lines, and under the $x-$axis, which is the comb part, we add many horizontal 
lines. Theorem 4.1 is saying  that in this case we still have  all the  main 
features of the HPHC walk.

\smallskip
\noindent
{\bf Example 3.} The easiest  example for which our Theorem  4.2  can be 
applied is the increasing gap case, namely when 
$$B=B(\alpha)=\{\pm [k^{\alpha}],k=0,1,2,...
\qquad {\rm with }\quad \alpha>1\}.$$
Then Theorem 4.2 holds true with $\beta=\frac{1}{\alpha}.$

\bigskip
It would be interesting to generalize Theorem 4.2  by replacing condition  
(\ref{beta}) with
\begg |B_n|\sim c n^{\beta}L(n)\endd
where $L(n)$ is a slowly varying function.

\bigskip

{\bf Acknowledgements} The authors are indebted to Mikl\'os Cs\"org\H{o} and 
P{\'al R\'ev\'esz 
for their 
\noindent
inspiration  and careful reading of our manuscript that greatly improved our 
presentation. We also would like to thank to our referee for his/her 
insightful suggestions which made  our presentation much nicer.

\end{document}